\documentclass[a4paper,12pt]{article}
\usepackage{amsmath}
\usepackage{latexsym}
\numberwithin{equation}{section}

\def\R{{\bf R}}

\def\N{{\bf N}}

\def\d{\displaystyle}
\def\e{{\varepsilon}}

\def\v#1{\mbox{\boldmath $#1$}}

\newtheorem{thm}{Theorem}[section]

\newtheorem{lem}{Lemma}[section]
\newtheorem{prop}{Proposition}[section]
\newtheorem{rem}{Remark}[section]

\title{The lifespan of solutions to wave equations with 
weighted nonlinear terms\\ in one space dimension}
\author{
Kyouhei Wakasa
\footnote{
Department of Mathematics, Hokkaido University, Sapporo, 060-0810, Japan. 
e-mail: wakasa@math.sci.hokudai.ac.jp.
}
}
\date{
\[
\begin{array}{l}
\mbox{\scriptsize{\bf Keywords:} nonlinear wave equation, lifespan, one space dimension}\\
\mbox{\scriptsize{\bf MSC2010:} Primary 35L71, 35E15, 
Secondary 35A01, 35A09, 35B44}\\
\end{array}
\]
}

\begin{document}
\maketitle
\begin{abstract}
In this paper, we consider the initial value problem for 
nonlinear wave equation with weighted nonlinear terms in one space dimension.
Kubo $\&$ Osaka $\&$ Yazici \cite{KOY13} studied global solvability of the 
problem under different conditions on the nonlinearity and initial data, 
together with an upper bound of the lifespan for the problem. 
The aim of this paper is to improve the upper bound of the lifespan and 
to derive its lower bound which shows the optimality of our new upper bound. 
\end{abstract}


\section{Introduction}
In this paper we consider the initial value problem for nonlinear wave equations:
\begin{equation}
\label{IVP}
\left\{
\begin{array}{l}
\d u_{tt}-u_{xx}=H(x,u(x,t)), \quad (x,t)\in\R \times[0,\infty),\\
u(x,0)=\e f(x),\ u_t(x,0)=\e g(x), \quad x\in\R,
\end{array}
\right.
\end{equation}
where $u=u(x,t)$ is a scalar unknown function of space-time variables, 
$(f,g)\in C^2(\R)\times C^1(\R)$ and $\e>0$ is a \lq\lq small" parameter. 
The nonlinear term, $H$ is given by 
\begin{equation}
\label{H}
H(x,u)=\frac{F(u(x,t))}{(1+|x|^2)^{(1+a)/2}},
\end{equation}
where $a\ge-1$ and $F(u)=|u|^{p}$ or $|u|^{p-1}u$ with $p>1$. 
Let us define the lifespan $T_\e$ of $C^2$-solution of (\ref{IVP}) by 
\[
\begin{array}{ll}
T_\e \equiv T_\e(f,g):=\sup\{T\in(0,\infty)\ :
&\mbox{There exists a unique solution}\\
&\mbox{\ $u\in C^2(\R\times[0,T))$ of (\ref{IVP})}\}
\end{array}
\]
with arbitrarily fixed $(f,g)$.
\par
First of all, we recall known results for the case $a=-1$ in general 
spatial dimensions:
\[
\left\{
\begin{array}{l}
u_{tt}-\Delta u=|u|^p\quad \mbox{\rm{in}}\quad \R^n\times[0,\infty),\\
u(x,0)=\e f(x),\ u_t(x,0)=\e g(x), \quad \mbox{\rm{for}}\ x\in\R^n,
\end{array}
\right.
\]
where $n\ge1$. 
When $n\ge2$, there exists a critical exponent $p_0(n)$ such that $T_\e=\infty$ 
for \lq\lq small"  $\e$ with compact support if $p>p_0(n)$, 
and $T_\e<\infty$ for \lq\lq positive" $(f,g)$ if $1<p\le p_0(n)$. 
Actually, $p_0(n)$ is a positive root of the quadratic equation $(n-1)p^2-(n+1)p-2=0$. 
See e.g. Introduction in Takamura $\&$ Wakasa \cite{TW14} for the details. 
\par
On the other hand, when $n=1$, $F(u)=|u|^p$, and $(f,g)$ has a compact support and satisfies 
some positivity assumption, Kato \cite{Kato80} showed that $T_\e<\infty$ 
for any $p>1$. 
The difference between the cases $n\ge2$ and $n=1$ 
comes from the fact that the solutions to the homogeneous wave 
equations has a decay estimate, $|u(x,t)|\le (t+1)^{-(n-1)/2}$. 
Especially, the solution does not have decay property when $n=1$.
\par
The result due to \cite{Kato80} motivates one to introduce a weight function 
$(1+x^2)^{-(1+a)/2}$ in the nonlinearity for getting a global solution.
Actually, Suzuki \cite{Suzuki10} showed that 
$T_\e=\infty$ with $F(u)=|u|^{p-1}u$ for $p>(1+\sqrt{5})/2$ and $pa>1$ 
if $f$ and $g$ are odd functions and $\e$ is small enough, and 
Kubo $\&$ Osaka $\&$ Yazici \cite{KOY13} have obtained the same conclusion 
for any $p>1$ satisfying $pa>1$.
On the other hand, they showed that $T_\e<\infty$ for $F(u)=|u|^p$ 
with $p>1$ and $a\ge-1$ if $(f,g)$ satisfies $f\equiv 0$, $g(x)\ge0$ for $x\in\R$, 
and $\int_{\delta/2}^{\delta}g(y)dy>0$ with some $0<\delta<1$. 
Also, they obtained an upper bound of the lifespan, $T_\e\le C\e^{-p^2}$, 
where $C$ is a positive constant independent of $\e$. 
However, this estimate is not sharp at least in the case of $a=-1$.
In fact, Zhou \cite{Z92_one} has obtained the following estimate of the lifespan $T_\e$ for any $p>1$,
\begin{equation}
\label{lifespan_est_one_dim.}
c\e^{-(p-1)/2}\le T_\e\le C\e^{-(p-1)/2}
\quad \mbox{if}\quad\d\int_{\R}g(x)dx\neq0,
\end{equation}
where $c$ and $C$ are positive constants independent of $\e$.
\par
Our purpose in this paper is to extend Zhou's result to the case where 
$a>-1$.
To obtain a blow-up result, we require the following assumptions on the data:
\begin{equation}
\label{blowup_asm}
\begin{array}{l}
\mbox{Let $f\equiv 0$ and $g\in C^{1}(\R)$ does not vanish identically.}\\ 
\mbox{Assume $g(x)\ge0$ for all $x\in\R$ and $\d \int_{-1}^{1}g(y)dy>0$.}
\end{array}
\end{equation}
Then, we have the following blow-up theorem.
\begin{thm}
\label{thm:upper_lifespan}
Let $a\ge-1$ and $F(u)=|u|^{p-1}u$ or $|u|^p$ with $p>1$. 
Assume (\ref{blowup_asm}). 
Then, there exist positive constants $\e_0=\e_0(g,a,p)$ and 
$C=C(g,a,p)$ such that 
\begin{equation}
\label{upper_lifespan_1}
T_\e\le 
\left\{
\begin{array}{ll}
C\e^{-(p-1)/(1-a)} & \mbox{if}\ -1\le a<0,\\
\phi^{-1}(C\e^{-(p-1)}) & \mbox{if}\ a=0,\\
C\e^{-(p-1)} & \mbox{if}\ a>0,
\end{array}
\right.
\end{equation}
holds for any $\e$ with $0<\e\le \e_0$, 
where $\phi=\phi(s)$ is a function defined by $\phi(s)=s\log(2+s)$ for $s\ge0$.
\end{thm}
The proof of this theorem done by an iteration argument concerning point-wise estimates. 
Such kind of framework was introduced by John \cite{J79} in three space dimensions. 
The first step of the iteration argument comes from the linear estimate 
of the solution to the homogeneous wave equation from below.
Kubo $\&$ Osaka $\&$ Yazici \cite{KOY13} obtained such an estimate only 
in a strip domain, $\{0\le x-t\le \delta/2\}$, where $0<\delta<1$ is a constant. 
On the other hand, we are able to show a similar estimate in unbounded domain, $\{t-x\ge1\}$. 
This improvement enable us to establish sharp upper bound of $T_\e$. 
See Lemma \ref{lem:lin_est} and Remark \ref{rem:lin_3dim} for details. 
\par
To show the optimality of the upper bounds in Theorem \ref{thm:upper_lifespan}, 
we require the following assumptions on $(f,g)$ 
\begin{equation}
\begin{array}{l}
\label{hypo_data}
\mbox{$f\in C^{2}(\R)$ and $g \in C^{1}(\R)$ satisfy $\|f\|_{L^{\infty}(\R)}<\infty$}\\
\mbox{and $\|g\|_{L^1(\R)}<\infty$.}
\end{array}
\end{equation}
Then, we have the following theorem.
\begin{thm}
\label{thm:lower_lifespan}
Let $a\ge-1$ and $F(u)=|u|^{p-1}u$ or $|u|^p$ with $p>1$. 
Assume (\ref{hypo_data}). 
Then, there exists a positive constant $c=c(f,g,a,p)$ such that 
\begin{equation}
\label{lower_lifespan_1}
T_\e\ge 
\left\{
\begin{array}{ll}
c\e^{-(p-1)/(1-a)} & \mbox{if}\ -1\le a<0,\\
\phi^{-1}(c\e^{-(p-1)}) & \mbox{if}\ a=0,\\
c\e^{-(p-1)} & \mbox{if}\ a>0,
\end{array}
\right.
\end{equation}
holds for $\e>0$, where $\phi$ is the function in Theorem \ref{thm:upper_lifespan}. 
\end{thm}

\begin{rem}
One can easily generalize the assumption on $F$ in Theorem \ref{thm:lower_lifespan} 
as follows:
\begin{equation}
\begin{array}{l}
\label{hypo_F}
\mbox{$F\in C^1(\R)$ satisfies $F(0)=F'(0)=0$ and}\\
\mbox{$|F'(s)|\le pA|s|^{p-1}$ for $s\in \R$, where $p>1$ and $A>0$.}
\end{array}
\end{equation}
\end{rem}
\par
This paper is organized as follows. In the next section, we prepare 
some notations. 
The upper bounds of the lifespan and lower bounds of the lifespan
are obtained in Section 3 and Section 4, respectively.

\section{Notations}
In this section, we give some notations and definitions. 
\par
We define
\begin{equation}
\label{lin_int}
u^0(x,t)=\frac{1}{2}\{f(x+t)+f(x-t)\}+\frac{1}{2}\int_{x-t}^{x+t}g(y)dy
\end{equation}
and
\begin{equation}
\label{duhamel}
L(V)(x,t)=\frac{1}{2}\iint_{D(x,t)}V(y,s)dyds
\end{equation}
for $V\in C(\R\times[0,\infty))$, 
where
\[
D(x,t)=\{(y,s)\in\R\times [0,\infty):0\le s\le t, x-t+s\le y \le x+t-s\}.
\]
For $(f,g)\in C^2(\R)\times C^1(\R)$, if $u\in C(\R\times[0,\infty))$ is a 
solution of 
\begin{equation}
\label{int_eq}
u(x,t)=\e u^0(x,t)+L(H(\cdot,u))(x,t), \quad (x,t)\in\R\times[0,\infty),
\end{equation}
then $u\in C^2(\R\times[0,\infty))$ is the solution to the initial value problem (\ref{IVP}).
\par
For $T>0$, we define the following domains:
\begin{equation}
\label{blowup_d}
\begin{array}{ll}
\Gamma_1=\{(x,t)\in[0,\infty)\times [0,T]: t-x\ge 1\},\\
\Gamma_2=\{(x,t)\in[0,\infty)\times [0,T]:x\ge t-x\ge 1\},\\
\Sigma_{j}=\{(x,t)\in[0,\infty)\times [0,T]:t-x\ge l_{j}\},
\end{array}
\end{equation}
where
\begin{equation}
\label{l_j}
\left\{
\begin{array}{l}
l_1=3\\
\d l_j=l_1+\sum_{k=1}^{j-1}2^{-(k-1)}=l_1+2\left(1-\frac{1}{2^{j-1}}\right)\quad \mbox{for}\ j\ge2. 
\end{array}
\right.
\end{equation}
\par
\section{Upper bound of the lifespan}
\par
In this section, we prove Theorem \ref{thm:upper_lifespan}. 
It is sufficient to show that the solution to the integral equation, 
\begin{equation}
\label{int_eq'}
u(x,t)=\e u^0(x,t)+\frac{1}{2}\iint_{D(x,t)}\frac{|u(y,s)|^pdyds}{(1+y^2)^{(1+a)/2}}, \quad (x,t)\in\R\times[0,\infty),
\end{equation}
blows up in finite time.
Because, if $u\in C(\R\times[0,\infty)$ is a solution of (\ref{int_eq'}), 
then $u$ satisfies $u(x,t)\ge0$ for $(x,t)\in\R\times[0,\infty)$ 
by the assumptions in (\ref{blowup_asm}). 
Therefore, this $u$ must solve the equation (\ref{int_eq}) 
with $F(u)=|u|^{p-1}u$ by the uniqueness of solutions to (\ref{IVP}). 
\par
Before proving Theorem \ref{thm:upper_lifespan}, we prepare the following lemmas:
\begin{lem}
\label{lem:seq}
Let $p>1$, $a\ge-1$ and let us define a sequence 
\begin{equation}
\label{C_j_1}
\left\{
\begin{array}{ll}
\d C_{a,j}=\exp\{p^{j-1}(\log(C_{a,1}F_{p,a}^{-S_j}E_{p,a}^{1/(p-1)}))-\log E_{p,a}^{1/(p-1)}\}\ (j\ge 2),\\
\d C_{a,1}=c_0^pk_{a}\e^p,
\end{array}
\right.
\end{equation}
where 
\begin{equation}
\label{E}
E_{p,a}
=\left\{
\begin{array}{ll}
(p-1)^2/(2^{a+5}p^2), & \mbox{if}\  -1\le a<0,\\
(p-1)^2/(2p^2), & \mbox{if}\ a=0,\\
(p-1)/(2^{a+2}p). & \mbox{if}\ a>0,
\end{array}
\right.
\end{equation}
\begin{equation}
\label{F}
F_{p,a}
=\left\{
\begin{array}{ll}
p^2, & \mbox{if}\  -1\le a\le 0,\\
2p & \mbox{if}\ a>0,
\end{array}
\right.
\end{equation}
\begin{equation}
\label{k_a}
k_{a}
=\left\{
\begin{array}{ll}
2^{-(a+4)}, & \mbox{if}\  -1\le a<0,\\
2^{-1}, & \mbox{if}\ a=0,\\
2^{-(a+2)}. & \mbox{if}\ a>0,
\end{array}
\right.
\end{equation}
and
\begin{equation}
\label{def_const_1}
S_j=\sum_{i=1}^{j-1}\frac{i}{p^i}. 
\end{equation}
Then, we have the following relation:
\begin{equation}
\label{C_j+1}
C_{a,j+1}=\frac{C_{a,j}^pE_{p,a}}{F_{p,a}^{j}}\quad (j\in \N).
\end{equation}
\end{lem}
\par\noindent
{\bf Proof.} 
First, we shall show (\ref{C_j+1}) for $j=1$. 
One can easily get
\[
\begin{array}{ll}
\d \log\left(\frac{C_{a,1}^pE_{p,a}}{F_{p,a}}\right)
=p\log(C_{a,1}F_{p,a}^{-1/p})+\log E_{p,a}\\
=p\log(C_{a,1}F_{p,a}^{-1/p}E_{p,a}^{1/(p-1)})-\log E_{p,a}^{1/(p-1)}
=\log C_{a,2}.
\end{array}
\]
Hence (\ref{C_j+1}) holds for $j=1$.
Next, we shall show (\ref{C_j+1}) for $j\ge2$.
Note that (\ref{C_j+1}) is equivalent to  
\[
\log C_{a,j+1}=p\log C_{a,j}-j\log F_{p,a} +\log E_{p,a}.
\]
By (\ref{C_j_1}) and the expression of $S_j$ in (\ref{def_const_1}), 
the right-hand side of this identity is equal to  
\[
\begin{array}{ll}
p^j\{\log (C_{a,1}F_{p,a}^{-S_j}E_{p,a}^{1/(p-1)})\}-p\log E_{p,a}^{1/(p-1)}
-j\log F_{p,a}+\log E_{p,a}\\
=p^j\{\log (C_{a,1}F_{p,a}^{-S_{j+1}}E_{p,a}^{1/(p-1)})\}
+p^j\log F_{p,a}^{j/p^j}-j\log F_{p,a}-\log E_{p,a}^{1/(p-1)}\\
=p^j\{\log (C_{a,1}F_{p,a}^{-S_{j+1}}E_{p,a}^{1/(p-1)})\}-\log E_{p,a}^{1/(p-1)}.\\
\end{array}
\]
Hence, we obtain (\ref{C_j+1}) by (\ref{C_j_1}) with $j$ replaced by $j+1$.
This completes the proof.
\hfill$\Box$
\par
Next, we derive a lower bound of the solution to (\ref{int_eq'}) which is 
a starting point of our iteration argument. 
\begin{lem}
\label{lem:lin_est}
Suppose that the assumptions in Theorem \ref{thm:upper_lifespan} are fulfilled. 
Let $u\in C(\R\times[0,T])$ be the solution of (\ref{int_eq'}). 
Then, $u$ satisfies
\begin{equation}
\label{linear_est}
u(x,t)\ge  \e c_0\quad \mbox{for}\quad (x,t)\in \Gamma_1,
\end{equation}
where $\d c_0=\frac{1}{2}\int_{-1}^{1}g(y)dy>0$ and $\Gamma_1(=\{(x,t)\in[0,\infty)\times [0,T]: t-x\ge 1\})$ is the
one in (\ref{blowup_d}).
\end{lem}
\par\noindent
{\bf Proof.} By (\ref{blowup_asm}) and (\ref{lin_int}), we get 
\[
\e u^0(x,t)=\frac{\e}{2}\int_{x-t}^{x+t}g(y)dy\ge \e c_0 \quad \mbox{for}\quad (x,t)\in \Gamma_1.
\]
Making use of the positivity of the second term of right-hand side in 
(\ref{int_eq'}), we have (\ref{linear_est}).
This completes the proof.
\hfill$\Box$
\begin{rem}
\label{rem:lin_3dim}
In three space dimensions, the following estimate which is necessary 
to get the first step of the iteration argument 
was obtained by John \cite{J79} in a strip domain: 
For $(x,t)\in S$, we have
\[
u^0(x,t)\ge Cr^{-1},
\]
where $r=|x|$, $C$ is a positive constant and 
$S=\{(r,t)\in(0,\infty)\times[0,\infty):\delta \le t-r\le \delta'\}$,
with some $\delta', \delta$ $(\delta'>\delta>0)$.
\par
On the contrary, our estimate holds in some domain without 
any restriction of upper bound for $t-x$. 
This is the key point to obtain sharp upper bound of $T_\e$.
\end{rem}

\par
Our iteration argument will be done by using the following estimates.
\begin{prop}
\label{prop:iteration_frame}
Suppose that the assumptions in Theorem \ref{thm:upper_lifespan} are fulfilled. 
Let $j\in\N$ and let $u\in C(\R\times[0,T])$ be the solution of (\ref{int_eq'}). 
Then, $u$ satisfies

\begin{equation}
\label{iteration_1}
u(x,t)\ge C_{a,j}\{(t-x)^{-(a+1)}(t-x-1)^2\}^{a_j}\quad \mbox{if}\ -1\le a <0,
\end{equation}
for $(x,t)\in \Gamma_{2}$, and 
\begin{equation}
\label{iteration_2}
u(x,t)\ge C_{0,j}\{(t-x-1)\log(1+x)\}^{a_j}\quad \mbox{if}\  a=0,
\end{equation}
for $(x,t)\in \Gamma_{1}$, and
\begin{equation}
\label{iteration_3}
u(x,t)\ge C_{a,j}(t-x-l_{j})^{a_j}\quad \mbox{if}\ a>0,
\end{equation}
for $(x,t)\in\Sigma_{j}$,
where $\Gamma_{1}$, $\Gamma_{2}$ and $\Sigma_{j}$ are defined in (\ref{blowup_d}). 
Here $C_{a,j}$ is the one in (\ref{C_j_1}) with $\d c_0=\frac{1}{2}\int_{-1}^{1}g(y)dy>0$ 
and $a_j$ is defined by 
\begin{equation}
\label{a_j}
a_j=\frac{p^j-1}{p-1}\quad (j\in\N).
\end{equation}

\end{prop}
\par\noindent
{\bf Proof.}
We shall show (\ref{iteration_1}), (\ref{iteration_2}) and (\ref{iteration_3}) by induction. 
Noticing that $u^0(x,t)\ge0$ for $(x,t)\in\R\times[0,\infty)$ 
and $(1+y^2)^{1/2}\le 1+|y|$, we get
\begin{equation}
\label{frame_1}
u(x,t)\ge \frac{1}{2}\iint_{D(x,t)}\frac{|u(y,s)|^p}{(1+|y|)^{1+a}}dyds
\quad \mbox{in}\ \R\times[0,\infty).
\end{equation}
\vskip10pt
\par\noindent
{\bf (i) Estimate in the case of $\v{-1\le a<0}$.}
\par
Let $(x,t)\in\Gamma_2$. Define 
\[
T_1(x,t):=\{(y,s)\in D(x,t): 1\le s-y\le t-x\le y, s+y\le t+x\}.
\]
Changing the variables in the integral of (\ref{frame_1}) by 
\begin{equation}
\label{alpha_beta}
\alpha=s+y,\  \beta=s-y
\end{equation}
and replacing the domain of integration by $T_1(x,t)$, we get 
\begin{equation}
\label{frame_1_1}
u(x,t)\ge 
\frac{1}{4}\int_{1}^{t-x}d\beta\int_{2(t-x)+\beta}^{t+x}
\frac{|u(y,s)|^p}{\left\{1+(\alpha-\beta)/2\right\}^{1+a}}d\alpha
\quad \mbox{in}\ \Gamma_2.
\end{equation} 
Making use of (\ref{linear_est}) and $T_1(x,t) \subset \Gamma_1$
for $(x,t)\in \Gamma_2$,
we have 
\[
u(x,t)\ge 
\frac{c_0^p\e^p}{4}\int_{1}^{t-x}d\beta\int_{2(t-x)+\beta}^{t+x}
\frac{d\alpha}{\left\{1+(\alpha-\beta)/2\right\}^{1+a}}
\quad \mbox{in}\ \Gamma_2.
\]
Note that $x\ge t-x$ is equivalent to $t+x\ge 3(t-x)$, we get 
\[
u(x,t)\ge \frac{c_0^p\e^p}{4}\int_{1}^{t-x}d\beta\int_{2(t-x)+\beta}^{3(t-x)}
\frac{d\alpha}{\left\{1+(\alpha-\beta)/2\right\}^{1+a}}
\quad \mbox{in}\ \Gamma_2.
\]
It follows from
\[
1+\frac{\alpha-\beta}{2}\le 1+\frac{3(t-x)-1}{2}\le2(t-x)
\]
for $\alpha\le 3(t-x)$, $\beta\ge 1$ and $t-x\ge1$
that 
\[
u(x,t)\ge \frac{c_0^p\e^p}{2^{a+3}(t-x)^{1+a}}
\int_{1}^{t-x}(t-x-\beta)d\beta=C_{a,1}\frac{(t-x-1)^2}{(t-x)^{1+a}}
\quad \mbox{in}\ \Gamma_2.
\]
Therefore, (\ref{iteration_1}) holds for $j=1$.
\par
Assume that (\ref{iteration_1}) holds. 
Noticing that $T_1(x,t) \subset \Gamma_2$ for $(x,t)\in \Gamma_2$
and putting (\ref{iteration_1}) into (\ref{frame_1_1}), we have
\[
u(x,t)
\ge\frac{C_{a,j}^p}{4}\int_{1}^{t-x}\frac{(\beta-1)^{2pa_j}}{\beta^{p(a+1)a_j}}d\beta
\int_{2(t-x)+\beta}^{t+x}\frac{d\alpha}{\left\{1+(\alpha-\beta)/2\right\}^{1+a}}
\quad \mbox{in}\ \Gamma_2.
\]
Analogously to the case of $j=1$, we get
\[
\begin{array}{lll}
u(x,t)
&\ge\d \frac{C_{a,j}^p}{2^{a+3}(t-x)^{(a+1)(pa_j+1)}}\int_{1}^{t-x}(\beta-1)^{2pa_j}d\beta
\int_{2(t-x)+\beta}^{3(t-x)}d\alpha\\
&=\d\frac{C_{a,j}^p}{2^{a+3}(t-x)^{(a+1)(pa_j+1)}}\int_{1}^{t-x}(\beta-1)^{2pa_j}(t-x-\beta)d\beta
\end{array}
\]
in $\Gamma_2$.
Making use of integration by parts to the integral above, we have 
\[
u(x,t)\ge\frac{C_{a,j}^p(t-x-1)^{2(pa_j+1)}}{2^{a+5}(pa_j+1)^2(t-x)^{(a+1)(pa_j+1)}}
\quad \mbox{in}\ \Gamma_2.
\]
Recalling the definition of $a_j$, we have 
\begin{equation}
\label{est_a_j}
a_{j+1}=pa_{j}+1\le \frac{p^{j+1}}{p-1}.
\end{equation}
Making use of (\ref{C_j+1}), we get 
\[
u(x,t)
\d\ge\frac{C_{a,j}^p(p-1)^2}{2^{a+5}p^{2(j+1)}}\cdot\frac{(t-x-1)^{2a_{j+1}}}{(t-x)^{(a+1)a_{j+1}}}
=C_{a,j+1}\frac{(t-x-1)^{2a_{j+1}}}{(t-x)^{(a+1)a_{j+1}}}
\]
in $\Gamma_2$. Therefore, (\ref{iteration_1}) holds for all $j\in\N$.
\vskip10pt
\par\noindent
{\bf (ii) Estimate in the case of $\v{a=0}$.}
\par
Let $(x,t)\in\Gamma_1$. Define 
\[
T_2(x,t):=\{(y,s)\in D(x,t): 1\le s-y\le t-x, s+y\le t+x, y\ge 0\}.
\]
Changing the variables by (\ref{alpha_beta}) in the integral of (\ref{frame_1}) 
and replacing the domain of integration by $T_2(x,t)$, we get 
\begin{equation}
\label{frame_1_2}
u(x,t)\ge\frac{1}{4}\int_{1}^{t-x}d\beta\int_{\beta}^{t+x}\frac{|u(y,s)|^p}{1+(\alpha-\beta)/2}d\alpha
\quad \mbox{in}\ \Gamma_1.
\end{equation}
By making use of (\ref{linear_est}) and $T_2(x,t) \subset \Gamma_1$ 
for $(x,t)\in \Gamma_1$, we get
\[
u(x,t)\ge\frac{c_0^p\e^p}{4}\int_{1}^{t-x}d\beta\int_{\beta}^{t+x}\frac{d\alpha}{1+(\alpha-\beta)/2}
\quad \mbox{in}\ \Gamma_1.
\]
Noticing that 
\[
\int_{\beta}^{t+x}\frac{d\alpha}{1+(\alpha-\beta)/2}=2\log\left(1+\frac{t+x-\beta}{2}\right)
\ge 2\log(1+x),
\]
for $\beta\le t-x$, 
we obtain
\[
u(x,t)\ge\frac{c_0^p\e^p}{2}\log(1+x)\int_{1}^{t-x}d\beta
=C_{0,1}(t-x-1)\log(1+x)\quad \mbox{in}\ \Gamma_1.
\]
Therefore, (\ref{iteration_2}) holds for $j=1$.
\par
Assume that (\ref{iteration_2}) holds. 
Noticing that $T_2(x,t) \subset \Gamma_1$ for $(x,t)\in \Gamma_1$
and putting (\ref{iteration_2}) into (\ref{frame_1_2}), we have
\[
u(x,t)\ge \frac{C_{0,j}^p}{4}\int_{1}^{t-x}(\beta-1)^{pa_j}d\beta
\int_{\beta}^{t+x}\frac{\left\{\log\left(1+(\alpha-\beta)/2\right)\right\}^{pa_j}d\alpha}{1+(\alpha-\beta)/2}
\quad \mbox{in}\ \Gamma_1.
\]
Analogously to the case of $j=1$, we get
\[
\begin{array}{llll}
u(x,t)&\d\ge \frac{C_{0,j}^p}{2(pa_j+1)}\int_{1}^{t-x}(\beta-1)^{pa_j}
\left\{\log\left(1+\frac{t+x-\beta}{2}\right)\right\}^{pa_j+1}d\beta&\\
&\d\ge \frac{C_{0,j}^p\left\{\log(1+x)\right\}^{pa_j+1}}{2(pa_j+1)}\int_{1}^{t-x}(\beta-1)^{pa_j}d\beta&
\end{array}
\]
in $\Gamma_1$. 
It follows from (\ref{est_a_j}) and (\ref{C_j+1}) that 
\[
\begin{array}{lll}
u(x,t)&\d\ge \frac{C_{0,j}^p(p-1)^2}{2p^{2(j+1)}}\cdot\left\{(t-x-1)\log(1+x)\right\}^{a_{j+1}}&\\
&=C_{0,j+1}\{(t-x-1)\log(1+x)\}^{a_{j+1}}&
\end{array}
\]
in $\Gamma_1$.
Therefore, (\ref{iteration_2}) holds for all $j\in\N$.

\vskip10pt
\par\noindent
{\bf (iii) Estimate in the case of $\v{a>0}$.}
\par
Let $(x,t)\in \Sigma_{1}$. Define 
\[
L_1(x,t):=\{(y,s)\in D(x,t): 1\le s-y\le t-x-2, 0\le y\le1\}.
\]
Changing the variables by (\ref{alpha_beta}) in the integral of (\ref{frame_1}) 
and replacing the domain of integration by $L_1(x,t)$, 
we get 
\[
u(x,t)\ge \frac{1}{4}\int_{1}^{t-x-2}d\beta
\int_{\beta}^{2+\beta}\frac{|u(y,s)|^pd\alpha}{\left\{1+(\alpha-\beta)/2\right\}^{1+a}}
\quad \mbox{in}\ \Sigma_{1}.
\]
By making use of (\ref{linear_est}) and $L_1(x,t) \subset \Gamma_1$ for $(x,t)\in \Sigma_1$, 
we have
\[
u(x,t)\ge \frac{c_0^p\e^p}{4}\int_{1}^{t-x-2}d\beta\int_{\beta}^{2+\beta}\d \frac{d\alpha}{\left\{1+(\alpha-\beta)/2\right\}^{1+a}}
\quad \mbox{in}\ \Sigma_{1}.
\]
It follows from $1+(\alpha-\beta)/2\le 2$ for $\alpha\le 2+\beta$ that
\[
u(x,t)\ge \frac{c_0^p\e^p}{2^{a+2}}
\int_{1}^{t-x-2}d\beta=C_{a,1}(t-x-3) \quad \mbox{in}\ \Sigma_{1}.
\]
Therefore, (\ref{iteration_3}) holds for $j=1$.
\par
Assume that (\ref{iteration_3}) holds. 
Let $(x,t)\in \Sigma_{j+1}$. Define 
\[
L_j(x,t):=\{(y,s)\in D(x,t): l_j\le s-y\le t-x-2^{-(j-1)},0\le y\le 2^{-j}\}
\] 
for $j\ge2$, where $l_j$ is defined in (\ref{l_j}). 
Making use of (\ref{alpha_beta}) and replacing the domain of integration 
in (\ref{frame_1}) by $L_j(x,t)$, we have
\[
u(x,t)\ge \frac{C_{a,j}^p}{4}\int_{l_j}^{t-x-2^{-(j-1)}}d\beta
\int_{\beta}^{2^{-(j-1)}+\beta}\frac{|u(y,s)|^pd\alpha}{\left\{1+(\alpha-\beta)/2\right\}^{1+a}}
\quad \mbox{in}\ \Sigma_{j+1}.
\]
Noticing that $L_j(x,t) \subset \Sigma_j$ for $(x,t)\in \Sigma_{j+1}$ and 
putting (\ref{iteration_3}) into the integral above, we have
\[
u(x,t)\ge \frac{C_{a,j}^p}{4}\int_{l_j}^{t-x-2^{-(j-1)}}(\beta-l_j)^{pa_j}d\beta
\int_{\beta}^{2^{-(j-1)}+\beta}\frac{d\alpha}{\left\{1+(\alpha-\beta)/2\right\}^{1+a}}
\]
in $\Sigma_{j+1}$.
Note that 
\[
1+\frac{\alpha-\beta}{2}\le 1+\frac{1}{2^j}\le 2
\]
for $\alpha\le 2^{-(j-1)}+\beta$, we get
\[
u(x,t)\ge \frac{C_{a,j}^p}{2^{a+2+j}}\int_{l_j}^{t-x-2^{-(j-1)}}(\beta-l_j)^{pa_j}d\beta
\quad \mbox{in}\ \Sigma_{j+1}.
\]
It follows from $\d l_j+2^{-(j-1)}=l_{j+1}$, (\ref{est_a_j}) and (\ref{C_j+1}) that 
\[
u(x,t)\ge \frac{(p-1)C_{a,j}^p}{2^{a+2+j}p^{(j+1)}}\cdot(t-x-l_{j+1})^{a_{j+1}}
=C_{a,j+1}(t-x-l_{j+1})^{a_{j+1}}
\]
in $\Sigma_{j+1}$. 
Therefore, (\ref{iteration_3}) holds for all $j\in\N$.
The proof of Proposition \ref{prop:iteration_frame} is now completed. 
\hfill $\Box$
\vskip10pt
\par\noindent
{\bf End of the proof of Theorem \ref{thm:upper_lifespan}.}\ 
Let $u\in C(\R\times[0,T])$ be the solution of the integral equation, (\ref{int_eq'}).
Setting $\d S=\lim_{j\rightarrow \infty}S_j$, we see from (\ref{def_const_1}) 
that $S_j\le S$ for all $j\in\N$. Therefore, (\ref{C_j_1}) yields 
\begin{equation}
\label{C_j_est}
\begin{array}{lll}
C_{a,j}&\ge\exp\{p^{j-1}\{\log(C_{a,1}F_{p,a}^{-S}E_{p,a}^{1/(p-1)})\}-\log E_{p,a}^{1/(p-1)}\}&\\
&=E_{p,a}^{-1/(p-1)}\exp\{p^{j-1}\{\log(C_{a,1}F_{p,a}^{-S}E_{p,a}^{1/(p-1)})\}\}.&
\end{array}
\end{equation}
\vskip10pt
\par\noindent
{\bf (i) The lifespan in the case of $\v{-1\le a<0}$.}
\par
We take $\e_0=\e_0(g,a,p)>0$ so small that 
\[
B_1\e_0^{-(p-1)/(1-a)}\ge4,
\]
where we set
\[
B_1=(c_0^p2^{-(a+4)+p(a-3)/(p-1)}p^{-2S}E_{p,a}^{1/(p-1)})^{-(p-1)/p(1-a)}>0. 
\]
Next, for a fixed $\e\in(0,\e_0]$, we suppose that $T$ satisfies 
\begin{equation}
\label{T_asm_1}
T>B_1\e^{-(p-1)/(1-a)}\ (\ge4).
\end{equation}
Combining (\ref{C_j_est}) with (\ref{iteration_1}), we have
\[
\begin{array}{ll}
\d u(x,t)\d\ge E_{p,a}^{-1/(p-1)}\exp\{p^{j-1}\{\log(C_{a,1}F_{p,a}^{-S}E_{p,a}^{1/(p-1)})\}\}\\
\d\qquad\times\left\{\frac{(t-x-1)^2}{(t-x)^{(1+a)}}\right\}^{(p^{j}-1)/(p-1)}
\end{array}
\]
in $\Gamma_2$.
Note that $\d t-x-1\ge (t-x)/2$ is equivalent to $t-x\ge2$. 
Furthermore, we have $(t/2,t)\in \Gamma_2$ for $t\in[4,T]$. 
Hence we get
\[
\begin{array}{llll}
\d u(t/2,t)
\!\!&\ge(2^{a-3}E_{p,a})^{-1/(p-1)}\exp\{p^{j-1}\{\log(2^{p(a-3)/(p-1)}C_{a,1}F_{p,a}^{-S}E_{p,a}^{1/(p-1)})\}\}\\
&\quad \times t^{(1-a)(p^{j}-1)/(p-1)}&\\
&=(2^{a-3}E_{p,a})^{-1/(p-1)}\exp\{p^{j-1}K_1(t)\}t^{-(1-a)/(p-1)}&
\end{array}
\]
for $t\in[4,T]$, where we set
\[
K_1(t)=\log \left(\e^pc_0^p2^{-(a+4)+p(a-3)/(p-1)}p^{-2S}E_{p,a}^{1/(p-1)}t^{p(1-a)/(p-1)}\right)
\]
(recall (\ref{F}) and (\ref{k_a})).
\par
By (\ref{T_asm_1}) and the definition of $B_1$, 
we have $K_1(T)>0$. 
Therefore we get $u(T/2,T)\rightarrow \infty$ as 
$j\rightarrow \infty$. 
Hence, (\ref{T_asm_1}) implies that $T_\e\le B_1\e^{-(p-1)/(1-a)}$ for $0<\e\le \e_0$. 
\vskip10pt
\par\noindent
{\bf (ii) The lifespan in the case of $\v{a=0}$.}
\par
We take $\e_1=\e_1(g,p)>0$ so small that 
\[
\phi^{-1}(B_2\e_1^{-(p-1)})\ge4,
\]
where $\phi$ is the one in Theorem \ref{thm:upper_lifespan} and 
\[
B_2=(c_0^p2^{-1-3p/(p-1)}p^{-2S}E_{p,0}^{1/(p-1)})^{-(p-1)/p}>0.
\]
Next, for a fixed $\e\in(0,\e_1]$, we suppose that $T$ satisfies 
\begin{equation}
\label{T_asm_2}
T>\phi^{-1}(B_2\e^{-(p-1)})\ (\ge4).
\end{equation}
Combining the estimates (\ref{C_j_est}) and (\ref{iteration_2}), we have
\[
\begin{array}{l}
\d u(t/2,t)
\d \ge (2^{-2}E_{p,0})^{-1/(p-1)}
\exp\{p^{j-1}\{\log(\e^pc_0^p2^{-1-2p/(p-1)}p^{-2S}E_{p,0}^{1/(p-1)})\}\}\\
\qquad\times\{t\log(1+t/2)\}^{(p^j-1)/(p-1)}
\end{array}
\]
for $4\le t \le T$.
Noticing that 
\[
\log\left(1+\frac{t}{2}\right)=\log(2+t)-\log2\ge\frac{\log(2+t)}{2}\quad \mbox{for}\ t\ge2,
\]
we get
\[
u(t/2,t)
\ge (2^{-3}E_{p,0})^{-1/(p-1)}\exp\{p^{j-1}K_2(t)\}\phi(t)^{-1/(p-1)}
\]
for $4\le t \le T$,
where we set
\[
\begin{array}{l}
\d K_2(t)=\log \left(\e^pc_0^p2^{-1-3p/(p-1)}p^{-2S}E_{p,0}^{1/(p-1)}\{\phi(t)\}^{p/(p-1)}\right).
\end{array}
\]
\par
Analogously to the case of $-1\le a<0$, 
we have $K_2(T)>0$ by (\ref{T_asm_2}) and the definition of $B_2$. 
Therefore we get $u(T/2,T)\rightarrow \infty$ as $j\rightarrow \infty$. 
Hence, (\ref{T_asm_2}) implies that $T_\e\le \phi^{-1}(B_2\e^{-(p-1)})$ for $0<\e\le \e_1$. 
\vskip10pt
\par\noindent
{\bf (iii) The lifespan in the case of $\v{a>0}$.}
\par
We take $\e_2=\e_2(g,a,p)>0$ so small that 
\[
B_3\e_2^{-(p-1)}\ge20,
\]
where we set
\[
B_3=(c_0^p2^{-(a+2)-2p/(p-1)}(2p)^{-S}E_{p,a}^{1/(p-1)})^{-(p-1)/p}>0. 
\]
Next, for a fixed $\e\in(0,\e_2]$, we suppose that $T$ satisfies 
\begin{equation}
\label{T_asm_3}
T>B_3\e^{-(p-1)}\ (\ge20).
\end{equation}
Combining the estimates (\ref{C_j_est}) with (\ref{iteration_3}), we have
\[
u(t/2,t)\ge (2^{-2}E_{p,a})^{-1/(p-1)}\exp\{p^{j-1}K_3(t)\}t^{-1/(p-1)}
\]
for $20\le t \le T$,
where we set
\[
K_3(t)=\log \left(\e^pc_0^p2^{-(a+2)-2p/(p-1)}(2p)^{-S}E_{p,a}^{1/(p-1)}t^{p/(p-1)}\right).
\]
\par
Since $K_3(T)>0$, by (\ref{T_asm_3}) and the definition of $B_3$, 
we get $u(T/2,T)\rightarrow \infty$ as $j\rightarrow \infty$. 
Hence, (\ref{T_asm_3}) implies that $T_\e\le B_3\e^{-(p-1)}$ for $0<\e\le \e_2$. 
Therefore, the proof of Theorem \ref{thm:upper_lifespan} is now completed.
\hfill$\Box$


\section{Lower bound of the lifespan}
In this section, we prove Theorem \ref{thm:lower_lifespan}. 
First of all, we introduce a Banach space 
\begin{equation}
X=\{u\in C(\R\times[0,T]): \|u\|_{L^{\infty}(\R\times[0,T])}<\infty\},
\end{equation}
which is equipped with a norm
\begin{equation}
\|u\|_{L^{\infty}(\R\times[0,T])}=\sup_{(x,t)\in\R\times [0,T]}|u(x,t)|.
\end{equation}
We shall construct a solution of the integral equation (\ref{int_eq}) in $X$ 
under suitable assumption on $T$ such as (\ref{det_lifespan_1}) below. 
Define a sequence of functions $\{u_n\}_{n\in\N}$ by 
\begin{equation}
\label{frame_exs}
u_n=u_0+L(H(\cdot,u_{n-1})), \quad u_0=\e u^0,
\end{equation}
where $L$, $H$ and $u^0$ are given by (\ref{duhamel}), (\ref{H}) and (\ref{lin_int}), 
respectively.
Since $\|u_0\|_{L^{\infty}(\R\times[0,T])}\le M\e$, where 
$M=\|f\|_{L^{\infty}(\R)}+\|g\|_{L^1(\R)}$ by (\ref{lin_int}), we have $u_0\in X$.
\par
The following {\it a priori} estimate plays a key role in the proof of 
Theorem \ref{thm:lower_lifespan}. 
\begin{lem}
\label{lem:frame_diff}
Let $V\in X$, $a\ge-1$, and let $D=D(\tau)$ is a function defined by 
\begin{equation}
\label{D}
D(\tau)=\left\{
\begin{array}{ll}
(1+\tau)^{1-a} & \mbox{if}\ -1\le a<0,\\
\phi(\tau) & \mbox{if}\ a=0,\\
\ 1+\tau & \mbox{if}\ a>0,
\end{array}
\right.
\end{equation}
for $\tau\ge0$, where $\phi$ is the one in Theorem \ref{thm:upper_lifespan}. 
Then, there exists a positive constant $C_a$ such that 
\begin{equation}
\label{diff_1}
\left\|L\left(\frac{V}{(1+|\cdot|^2)^{(1+a)/2}}\right)\right\|_{L^{\infty}(\R\times[0,T])}
\le C_a D(T)\|V\|_{L^{\infty}(\R\times[0,T])}.
\end{equation}
\end{lem}
{\bf Proof.} Noticing that $(1+y^2)\ge (1+|y|)^2/2$, 
the left-hand side in (\ref{diff_1}) is dominated by 
\[
C_a\|V\|_{L^{\infty}(\R\times[0,T])}\iint_{D(x,t)}\frac{dyds}{\langle y\rangle^{1+a}},
\]
where we set $\langle y\rangle=1+|y|$.
Thus, it is enough to show the inequality,
\begin{equation}
\label{basic_est}
I(x,t)\le C_aD(T)\quad \mbox{for}\ (x,t)\in \R\times[0,T],
\end{equation}
where we set
\[  
I(x,t)=\iint_{D(x,t)}\frac{dyds}{\langle y\rangle^{1+a}}.
\]
\par
We may assume $x\ge0$.
Because $I(x,t)$ is an even function with respect to $x$. 
When $t\ge x\ge0$, we divide the integral domain $D(x,t)$ into two parts 
$D_j(x,t)$ $(j=1,2)$, where
\[
\begin{array}{ll}
D_1(x,t)=\{(y,s)\in\R\times [0,\infty): 0\le s\le t-x,x-t+s\le y\le t-x-s\},\\
D_2(x,t)=\{(y,s)\in[0,\infty)^2: 0\le s\le t,|x-t+s|\le y\le x+t-s\}.
\end{array}
\] 
Namely, we set 
\[
I_j(x,t)=\iint_{D_j(x,t)}\frac{1}{\langle y\rangle^{1+a}}dyds\quad (j=1,2),
\]
so that $\d I(x,t)=I_1(x,t)+I_2(x,t)$. 
We shall estimate $I_1$. 
Since $\langle y\rangle$ is an even function, 
we obtain
\[
I_1(x,t)=2\int_{0}^{t-x}ds\int_{0}^{t-x-s}\frac{dy}{(1+y)^{1+a}}
\quad \mbox{for}\ t\ge x\ge0.
\]
Then, the $y$-integral is dominated by
\[
\left\{
\begin{array}{ll}
-a^{-1}(1+t-x)^{-a} & \mbox{if}\ a<0,\\
\log(1+t-x) & \mbox{if}\ a=0,\\
a^{-1} & \mbox{if}\ a>0.
\end{array}
\right.
\]
Hence, we get
\[
I_1(x,t)\le C_aD(t-x)\le C_aD(T)\quad \mbox{for}\ 0\le x \le t\le T.
\]
Next, we shall estimate $I_2$. It follows that
\[
I_2(x,t)=\int_{0}^{t}ds\int_{|x-t+s|}^{t-x-s}\frac{dy}{(1+y)^{1+a}}
\le \int_{0}^{t}ds\int_{0}^{t+x-s}\frac{dy}{(1+y)^{1+a}}
\]
for $t\ge x\ge0$, and that the $y$-integral is dominated by
\[
\left\{
\begin{array}{ll}
-a^{-1}(1+t+x)^{-a} & \mbox{if}\ a<0,\\
\log(1+t+x) & \mbox{if}\ a=0,\\
a^{-1} & \mbox{if}\ a>0.
\end{array}
\right.
\]
Noticing that 
\[
\log(1+2t)\le \log2+\log(2+t)\le 2\log (2+t)\quad \mbox{for}\ t\ge0, 
\]
we get 
\[
I_2(x,t)\le C_aD(t+x)\le C_aD(T)\quad \mbox{for}\ 0\le x \le t\le T.
\]
\par
When $x\ge t$, we have 
\[
\begin{array}{ll}
\d I(x,t)\le\int_{0}^{t}\frac{ds}{(1+s)^{1+a}}\int_{x-t+s}^{x+t-s}dy
\le 2t\int_{0}^{t}\frac{ds}{(1+s)^{1+a}}\le C_aD(T).
\end{array}
\]
Therefore, the proof of Lemma \ref{lem:frame_diff} is ended.
\hfill$\Box$
\vskip10pt
\par
Now, we move on to the proof of Theorem \ref{thm:lower_lifespan}. 
First of all, we take $T>0$ such that 
\begin{equation}
\label{det_lifespan_1}
2^{p+1}pC_aD(T)M^{p-1}\e^{p-1}\le 1,
\end{equation}
where $C_a$ is the one in Lemma \ref{lem:frame_diff}.
We shall show
\begin{equation} 
\label{bound_u}
\|u_n\|_{L^{\infty}(\R\times[0,T])}\le 2M\e\quad (n\in \N),
\end{equation}
by induction. Assume that $\|u_{n-1}\|_{L^{\infty}(\R\times[0,T])}\le 2M\e\ (n\ge2)$.
It follows from (\ref{frame_exs}) and Lemma \ref{lem:frame_diff} that 
\[
\begin{array}{lll}
\|u_n\|_{L^{\infty}(\R\times[0,T])}&\le \|u_0\|_{L^{\infty}(\R\times[0,T])}
+\|L(H(\cdot,u_{n-1}))\|_{L^{\infty}(\R\times[0,T])}&\\
&\le M\e+C_aD(T)\|u_{n-1}\|_{L^{\infty}(\R\times[0,T])}^p.&
\end{array}
\]
The assumption of the induction yields that
\[
\|u_n\|_{L^{\infty}(\R\times[0,T])}\le M\e+C_a(2M\e)^pD(T).
\]
This inequality shows (\ref{bound_u}), provided (\ref{det_lifespan_1}) holds.
\par
Next we shall estimate the differences of $\{u_n\}_{n\in\N}$.
Since
\[
\begin{array}{lll}
|H(y,u_{n})-H(y,u_{n-1})|
\!\!&\d\le \frac{p}{(1+y^2)^{(1+a)/2}}(|u_{n-1}(y,s)|^{p-1}+|u_{n}(y,s)|^{p-1})&\\
&\quad\times|u_{n}(y,s)-u_{n-1}(y,s)|&
\end{array}
\]
for $(y,s)\in\R\times[0,\infty)$, we see from Lemma \ref{lem:frame_diff} that
\[
\begin{array}{l}
\|u_{n+1}-u_{n}\|_{L^{\infty}(\R\times[0,T])}\\
\le pC_aD(T)(\|u_n\|_{L^{\infty}(\R\times[0,T])}^{p-1}+\|u_{n-1}\|_{L^{\infty}(\R\times[0,T])}^{p-1})\|u_n-u_{n-1}\|_{L^{\infty}(\R\times[0,T])}.
\end{array}
\]
Making use of (\ref{bound_u}), we have 
\[
\|u_{n+1}-u_{n}\|_{L^{\infty}(\R\times[0,T])} \le \frac{1}{2} \|u_{n}-u_{n-1}\|_{L^{\infty}(\R\times[0,T])}\quad\mbox{for}\ n\in\N
\]
provided (\ref{det_lifespan_1}) holds.
Hence, we obtain 
\[
\|u_{n+1}-u_n\|_{L^{\infty}(\R\times[0,T])}
\le \frac{1}{2^{n}}\|u_1-u_0\|_{L^{\infty}(\R\times[0,T])}\quad\mbox{for}\ n\in \N.
\]
Therefore, $\{u_n\}_{n\in\N}$ is a Cauchy sequence in $X$ provided (\ref{det_lifespan_1}) holds.
Since $X$ is complete, there exists $u\in X$ such that $u_n$ converges uniformly 
to $u$ in $X$. 
Therefore, by taking limits under the integral sign, 
$u$ satisfies the integral equation (\ref{int_eq}), 
so that $u$ is the $C^2$-solution of (\ref{IVP}).
Hence, the proof of Theorem \ref{thm:lower_lifespan} is completed.
\hfill$\Box$

\section*{Acknowledgment}
The author would like to express his gratitude to Professor Hideo Kubo for 
his suggestion of this problem and valuable advice. It is also a pleasure 
to acknowledge Professor Hiroyuki Takamura for his useful comments. 
The author is supported by Grant-in-Aid for Scientific Research 
of JSPS Fellow No. 26-2330.


\bibliographystyle{plain}

\end{document}